\documentclass[12pt,a4paper]{article}
\usepackage[cp1251]{inputenc}
\usepackage[psamsfonts]{amssymb}
\usepackage{amscd,amssymb,hyperref}
\usepackage{color}   

\usepackage{ulem}

\begin{document}
\thispagestyle{empty}

\title{On the lower central  series of Baumslag-Solitar groups }

\author{{V. G. Bardakov, M. V. Neshchadim}}%

\maketitle {\small}

\begin{quote}
\noindent{\sc Abstract.} 
We find the lower central series for  residually nilpotent Baumslag-Solitar groups,  and
find the intersection of all terms of the lower central series. Also, we find non-abelian Bauslag-Solitar groups for which the lower central series  has  length 2. 
For some Baumslag-Solitar groups  a connection  is found between the intersection of all  subgroups of finite index and the intersection of all terms of the lower central series. 

Keywords: Baumslag-Solitar groups, residually nilpotent group, lower central series, residually finite group. 
\end{quote}

\medskip

\begin{center}
{\bf Introduction}
\end{center}

Bauslag-Solitar  group (BSG)  may be presented in the form
$$
BS(m,n)=\langle \, a, t \, \| \, t^{-1} a^m t = a^n \, \rangle,~m, n \in \mathbb{Z} \setminus \{ 0 \}.
$$
These groups were introduced in  \cite{BS} as a set of groups which contains  non-hopfian groups. Recall that a group $G$ is called {\it non-hopfian} if it contains a proper quotient that is isomorphic to $G$. For example,  $BS(2, 3)$ is non-hopfian.  Since non-hopfian groups are not residually finite, then in Baumslg-Solitar groups there exist non-residually finite groups with one defining relation.
A survey of residual properties of BSG can be found in  \cite{Mol, Mol3}.

In the present paper we  study  the lower central series of BSG.
In particular, we describe the lower central series of residually nilpotent BSG.
For non-residually nilpotent groups it is interesting to find the lower central series to understand the reason of  non-residual nilpotence.
The reasons   may be one of two: the lower central series  stabilizes at some finite step, or all its terms are different, but their intersection is non-trivial. In the last case, it is interesting to study the transfinite lower central series.  We construct non-solvable BSG which have  lower central series of length 2. 

D.~I.~Moldavanskii \cite{Mol4} proved that the Baumslag-Solitar group  $BS(m, n)$, $0 < m \leq |n|$, is residually nilpotent if and only if 
 $m =1$, $n\not=2$ or $n = \varepsilon m$ and  $m > 1$ is a power of some prime number, 
$\varepsilon = \pm 1$. We find the lower central series for any group from this class.
Also, we prove that a Baumslag-Solitar group is torsion-free residually nilpotent  if and only if $n = 1$; that is, if and only if $BS(m, n)$ is isomorphic to  $\mathbb{Z} \times \mathbb{Z}$.

 It is obvious that the quotient of any group by the intersection of all terms of the lower central series is the maximal  residually nilpotent  quotient. 
We find the intersection of all terms of the lower central series for some Baumslag-Solitar groups and find some connection between intersection of all terms of the lower central series and the intersection of all subgroups of finite index.

{\bf Asknolegement.} We  thank F.~A.~Dudkin for the survey talk on  Baumslag-Solitar groups which he  gave at the seminar  ``Evariste Galois'' in Novosibirsk State University. Also, we thank D.~I.~Moldavanskii, who sent  us the unpublished paper 
 \cite{Mol4}  and read the first version of the present paper, and corrected some misprints. 
We thank D.~N.~Azarov for  interesting discussions. We thank I. Marshall for language checking.

\begin{center}
{\bf \S~1. Preliminaries }
\end{center}

The following definitions  can be found in the paper of A.~I.~Malcev  \cite{Mal}.
Let $\mathfrak{C}$ be some class of groups. A group   $G$ is said to be {\it  residually  $\mathfrak{C}$-group}
 or simply  $\mathfrak{C}$--{\it residual}, if for any non-identity  element 
 $g \in G$ there exists a homomorphism  $\varphi$ of  $G$ to some group from  $\mathfrak{C}$ such that  $\varphi(g) \not= 1$. If
$\mathfrak{C}$ is the class of all finite groups, then  $G$ is called  {\it residually finite}
(abbreviated  RF).  If $\mathfrak{C}$ is the class of finite  $p$-groups, then $G$ is said to be $p$-{\it residually finite}
(abbreviated  RF$_p$). If $\mathfrak{C}$ is the class of nilpotent groups, then   $G$ is  said to be 
 {\it residually nilpotent} (abbreviated RN). If $\mathfrak{C}$ is the class of torsion-free nilpotent groups, then   $G$ is  said to be  
 {\it torsion free residually nilpotent}.
We shall denote by  $r \mathfrak{F}$ the class of RF groups,
by $r \mathfrak{F}_p$, the class of RF$_p$ group, by $\{r \mathfrak{F}_p \}_p$
the class of groups which are RF$_p$ for some prime $p$,
by $r \mathfrak{N}$, the class of RN groups. It is clear that for finitely generated groups the following inclusions hold:
$$
r \mathfrak{F}_p \subset \{r \mathfrak{F}_p \}_p \subset r \mathfrak{N} \subset r \mathfrak{F}.
$$
It is easy to see,  that all inclusions are strict. Note that in general  the inclusion  $r \mathfrak{N} \subseteq r \mathfrak{F}$ is not true.
Indeed, the quasi-cyclic group   $C_{p^{\infty}}$ is not residually finite but is abelian and hence is residually nilpotent.

Let  $A, B \subseteq G$ be two subsets of a group $G$, then by  $[A, B]$ we shall denote the subgroup which is generated by all commutators 
 $[a, b] = a^{-1} b^{-1} a b$, $a \in A$, $b \in B$.
 The sequence of groups 
$$
G = \gamma_1 G \geq \gamma_2 G \geq \ldots,~~~ \mbox{where}~~\gamma_{i+1} G = [\gamma_i G, G], ~~i = 1, 2, \ldots ,
$$
 is called the \textit{central series} of $G$.
In particular, $G' = \gamma_2 G = [G, G]$ is the commutator subgroup of  $G$. The group $G$ is residually nilpotent if and only if the intersection of all terms of the lower central series is trivial;  that is,
$$
\gamma_{\omega} G = \bigcap_{i=1}^{\infty} \gamma_i G = 1.
$$
The {\it length of the lower central series} of  $G$ is the minimal number  $N$  for which  $\gamma_N G = \gamma_{N+1} G$. 
The length is equal to some natural number or to some  limiting ordinal  $\omega$.

Since the groups  $BS(m,n)$, $BS(n,m)$, $BS(-m,-n)$ are pairwise isomorphic,  we shall assume  for the rest of the paper  without loss of generality that the parameters  satisfy the inequalities   $0 < m \leq |n|$.

 Necessary and sufficient conditions  for residual finiteness  of  $BS(m,n)$ were found in the paper of Baumslag-Solitar  \cite{BS} and Meskin \cite{Mes}.
 A simpler proof  can be find in  \cite{Mol}.

\medskip

PROPOSITION A (\cite{BS}, \cite{Mes}).
{\it  $BS(m,n)$ is residually finite if and only if  $m = 1$ or $|n| = m > 1$.}

\medskip

 Necessary and sufficient conditions  for  $p$-residual finiteness of  $BS(m,n)$  were found by D.~I.~Moldavanskii \cite{Mol1} (see also \cite{Mol}).

\medskip

PROPOSITION B (\cite{Mol1, Mol}).
{\it  $BS(m,n)$ is $p$-residually finite if and only if   $m = 1$ and  $n \equiv 1~(mod~p)$, or $n = m = p^r$ for some $r \geq 0$, or $n = -m$ for $m=2^r$, $r \geq 0$.}

\medskip

Since every $p$-residually finite group is residually nilpotent, Proposition B implies that the following groups are residually nilpotent: 
$$
BS(1, 1 + pk),~~k \in \mathbb{Z},~~BS(p^r, \pm p^r),~r \geq 0.
$$
We see that this list  is different from the list of residually finite groups. It does not contain the groups  $BS(m, m)$, where $m > 1$ and is not a power of some prime number. 

\medskip

\begin{center}
{\bf \S~2. Lower central series}
\end{center}

{\bf  Solvable Baumslag-Solitar groups.}  Any solvable Baumslag-Solitar  group has the form 
$$
BS(1,n)=\langle \, a, t \, \| \, t^{-1} a t=a^n \, \rangle .
$$
It is easy to show (see, for example, \cite{Ba}), that any element  $g$ of  $BS(1,n)$ can be  presented in a unique manner as a product 
$$
g= t^{k}a^{l}t^{-r}, \quad k,r \in \mathbb{N}\cup \left\{ 0 \right\}, l \in \mathbb{Z},
$$
where $n \not|~ l$ for $r > 0$.

For example,
$$
t a t a t^{-2 }= t^2 (t^{-1} a t ) a t^{-2} = t^2 a^{n+1} t^{-2}.
$$

It is evident that for  $n = 1$ the group $BS(1,1)$ is the free abelian of rank  2 and hence is torsion-free residually nilpotent. 
If $n = -1$, then  $BS(1,-1)$ is the fundamental group of  the Klein bottle, for which  the following proposition  is well-known and easy to check.

\medskip

PROPOSITION 1.
{\it The group
$$
G=BS(1,-1)=\left\langle \, a, t \, \| \, t^{-1}at=a^{-1} \, \right\rangle
$$
is residually nilpotent and  the terms of  its lower central series have the form:}
$$
\gamma_{s+1} G= \left\langle \, a^{2^s} \, \right\rangle,\quad s \geq 1,
$$
$$
\gamma_{\omega} G= \bigcap\limits_{s \geq 1} \gamma_s G =1.
$$

\medskip

The following theorem describes the structure of the lower central series of  solvable Baumslag-Solitar groups and generalises the last proposition.

\medskip

THEOREM 1. {\it Let $G = BS(1,n)$, $n \not= 1$. Then}

1) $\gamma_i G$, $i > 1$, {\it consists of elements}
$$
t^{l} a^{\alpha (n-1)^{i-1}} t^{-l}, \quad l \in \mathbb{N}\cup \left\{ 0 \right\},~~\alpha \in \mathbb{Z};
$$

2) {\it  $\gamma_i G$, $i > 1$, is isomorphic to the additive group  $(n-1)^{i-2}\mathbb{Z}[1 / n]$.
In particular,  $G$ is a metabelian group;}

3) $G = G' \leftthreetimes \mathbb{Z}$,
{\it where $\mathbb{Z}$ is generated by  the element  $t$;
the quotient  $\gamma_i G / \gamma_{i+1} G$ is isomorphic to $\mathbb{Z}_{n-1}$ for all  $i > 1$.}

\medskip

PROOF.
1) Since 
$$
a^{n-1} = a^{-1} t^{-1} a t = [a, t],
$$
the subgroup  $H$ which is generated by elements  $t^{l}a^{n-1}t^{-l}$,
lies in the commutator subgroup  $G'$. Let us show that any element of  $H$ has the form
$$
t^{l}a^{\alpha (n-1)}t^{-l},~~ l \in \mathbb{N}\cup \left\{ 0 \right\},~~\alpha \in \mathbb{Z}.
$$
Indeed, take two elements of $H$:
$$
h_1 = t^{l_1}a^{\alpha_1 (n-1)}t^{-l_1},~~ h_2 = t^{l_2}a^{\alpha_2 (n-1)}t^{-l_2}
$$
and find their product. If $l_1 < l_2$, then
$$
h_1 \cdot h_2 = t^{l_1}a^{\alpha_1 (n-1)}t^{-l_1} \cdot t^{l_2}a^{\alpha_2 (n-1)}t^{-l_2} =
t^{l_1}a^{\alpha_1 (n-1)}t^{l_2-l_1}a^{\alpha_2 (n-1)}t^{-l_2} =
$$
$$
= t^{l_2}a^{\alpha_1 (n-1) n^{l_2 - l_1} + \alpha_2 (n-1)}t^{-l_2} =
t^{l_2}a^{(n-1) (\alpha_1 n^{l_2 - l_1} + \alpha_2)}t^{-l_2}.
$$
If $l_1 > l_2$, then
$$
h_1 \cdot h_2 = t^{l_1}a^{\alpha_1 (n-1)}t^{-l_1} \cdot t^{l_2}a^{\alpha_2 (n-1)}t^{-l_2} =
t^{l_1}a^{\alpha_1 (n-1)}t^{-(l_1-l_2)}a^{\alpha_2 (n-1)}t^{-l_2} =
$$
$$
= t^{l_1} a^{\alpha_1 (n-1) + \alpha_2 (n-1) n^{l_1 - l_2}} t^{-l_1} =
t^{l_1}a^{(n-1) (\alpha_1 + \alpha_2 n^{l_1 - l_2})}t^{-l_1}.
$$
If $l_1= l_2$, then the result is evident. 

Let us show that  $H$ is normal in  $G$. Indeed,
$$
t^{-1} \left( t^{l}a^{n-1}t^{-l} \right) t = t^{l}a^{n(n-1)}t^{-l}=(t^{l}a^{n-1}t^{-l})^n,
$$
$$
t \left( t^{l}a^{n-1}t^{-l} \right)t^{-1} = t^{l+1}a^{n(n-1)}t^{-l-1},
$$
$$
a^{-1} \left( t^{l}a^{n-1}t^{-l} \right)a =
t^{l} (t^{-l}a^{-1}t^{l}) a^{n-1} (t^{-l}a^{-1}t^{l}) t^{-l}=
t^{l}a^{-n^l}a^{n-1} a^{n^l}t^{-l}=t^{l}a^{n-1}t^{-l}.
$$

From the last relation: 
$[a, t^{l}a^{n-1}t^{-l}]=1$
follows that  $H$ is an abelian subgroup.
Since
$$
\langle \, t, a \, \| \, t^{-1}at=a^n, \,\, a^{n-1}=1 \, \rangle =
\langle \, t, a \, \| \, t^{-1}at=a, \,\, a^{n-1}=1 \, \rangle,
$$
the quotient  $G/H$ is abelian and hence  $G' \leq H$.

Further, using the induction by  $i\geq 2$.
Suppose that we have proven that 
$$
\gamma_i G=
\left\{\, t^{l} a^{\alpha (n-1)^{i-1}} t^{-l}\, | \, l \in \mathbb{N}\cup \left\{ 0 \right\},
\,\,\alpha \in \mathbb{Z} \, \right\}.
$$
Since  $[G',a]=1$ and  $G''=1$, then $\gamma_{i+1} G\subseteq [\gamma_i G,t]$.
We have
$$
[t^{l} a^{\alpha (n-1)^{i-1}} t^{-l},t]=
t^{l} a^{-\alpha (n-1)^{i-1}} t^{-l} t^{-1} t^{l} a^{\alpha (n-1)^{i-1}} t^{-l}t=
t^{l} a^{\alpha (n-1)^{i}} t^{-l}.
$$
As above, we can check that the product of elements of the form 
$t^{l} a^{\alpha (n-1)^{i}} t^{-l}$, $l \in \mathbb{N}\cup \left\{ 0 \right\}$,
$\alpha \in \mathbb{Z}$ has the same form. Hence, 
$$
\gamma_{i+1} G=
\left\{\, t^{l} a^{\alpha (n-1)^{i}} t^{-l}\, | \, l \in \mathbb{N}\cup \left\{ 0 \right\},
\,\,\alpha \in \mathbb{Z} \, \right\}
$$
and the claim 1) is proved.

2) Let us show that  
$G'\cong \mathbb{Z}[1 / n]$.
To this end, denote  $b=a^{n-1}$ and define the map 
$$
\varphi : G'\longrightarrow \mathbb{Z}\left[ 1 / n \right]
$$
by the rule 
$$
t^{l} b^{\alpha} t^{-l} \longmapsto \frac{\alpha}{n^l}, \quad l \in \mathbb{N}\cup \left\{ 0 \right\}.
$$
We see that this map is onto.
Since, in  $G'$ any element  can be uniquely presented in the form 
$t^{l} b^{\alpha} t^{-l}$, $l \in \mathbb{N}\cup \left\{ 0 \right\}$, $\alpha \in \mathbb{Z}$,
where $n \not| \alpha$ for  $l > 0$, we  find that the map $\varphi$ is an embedding.
Let us prove that  $\varphi$ preserves operations.
As above, take two elements in   $G'$:
$$
h_1 = t^{l_1}a^{\alpha_1 (n-1)}t^{-l_1},~~ h_2 = t^{l_2}a^{\alpha_2 (n-1)}t^{-l_2}
$$
and using the above calculations, we get 
$$
\varphi (h_1 \cdot h_2) = \frac{\alpha_1 n^{l_2 - l_1} + \alpha_2}{n^{l_2}}~~\mbox{for}~~l_1 < l_2;
$$
$$
\varphi (h_1 \cdot h_2) = \frac{\alpha_1 + \alpha_2 n^{l_1 - l_2}}{n^{l_1}}~~\mbox{for}~~l_1 > l_2.
$$
On the other side,
$$
\varphi (h_1) \cdot \varphi (h_2) = \frac{\alpha_1}{n^{l_1}} + \frac{\alpha_2}{n^{l_2}} =
\frac{\alpha_1 n^{l_2} + \alpha_2 n^{l_1}}{n^{l_1+l_2}},
$$
that is equal to 
$$
\frac{\alpha_1 n^{l_2-l_1}+\alpha_2}{n^{l_2}} ~~\mbox{for}~~l_1 < l_2;
$$
$$
\frac{\alpha_1 + \alpha_2 n^{l_1 - l_2}}{n^{l_1}} ~~\mbox{for}~~l_1 > l_2.
$$
 Comparing these terms, we see that  $\varphi$ is a homomorphism.
Hence, $\varphi$ is an isomorphism.

Further, acting by $\varphi$
on $\gamma_i G$, we get the  induced  isomorphism
$\gamma_i G \cong (n-1)^{i-2}\mathbb{Z}[1 / n]$.
The claim  2) is proved.

3)
 We need to prove the isomorphism
$\gamma_i G / \gamma_{i+1} G\cong \mathbb{Z}_{n-1}$ for all  $i > 1$.

Let $x_i$, $i \in \mathbb{N}\cup \left\{ 0 \right\}$
be free generators of the  free abelian group. Then the group  $\mathbb{Z}[1 / n]$ has the presentation
$$
\mathbb{Z}\left[1 / n \right]\cong
\left\langle \, x_i,\,\, i \in \mathbb{N}\cup \left\{ 0 \right\} \,
\| \, nx_i=x_{i-1},~~ i \in \mathbb{N}, \,\, [x_k, x_l] = 1 ~\mbox{for all}~k, l \, \right\rangle,
$$
where the isomorphism is defined by the map  $n^{-i} \to x_i$, $i \in \mathbb{N}\cup \left\{ 0 \right\}$.
Hence,
$$
\mathbb{Z}\left[1 / n\right] / (n-1) \mathbb{Z}\left[1 / n\right] \cong
\left\langle \, x_i,\,\, i \in \mathbb{N}\cup \left\{ 0 \right\} \,
\| \, (n-1)x_0=0,\,\, nx_i=x_{i-1},\,\, (n-1)x_i=0,\right.
$$
$$
\left. i \in \mathbb{N}, \,\, [x_k, x_l] = 1 ~\mbox{for all}~k, l \, \right\rangle\cong
$$
$$
\cong
\left\langle \, x_i,\, i \in \mathbb{N}\cup \left\{ 0 \right\} \,
\| \, (n-1)x_0=0,\, nx_i=x_{i-1},\, nx_i=x_i,\, i \in \mathbb{N}, \, [x_k, x_l] = 1 ~\mbox{for all}~k, l \, \right\rangle\cong
$$
$$
\cong
\left\langle \, x_0 \,
\| \, (n-1)x_0=0 \, \right\rangle\cong \mathbb{Z}_{n-1}.
$$
In particular, we have the equality
$$
\mathbb{Z}\left[1 / n\right] =
\left\langle x_0,\,\, (n-1) \mathbb{Z}\left[1 / n\right]\right\rangle.
$$
Hence, after  multiplying by   $(n-1)^i$, we get 
$$
(n-1)^i\mathbb{Z}\left[1 / n\right] =
\left\langle (n-1)^ix_0,\,\, (n-1)^{i+1} \mathbb{Z}\left[1 / n\right]\right\rangle.
$$
If for some  $1\leq k\leq n-1$ the following inclusion holds 
$$
k(n-1)^ix_0\in (n-1)^{i+1} \mathbb{Z}\left[1 / n\right],
$$
then
$$
kx_0\in (n-1) \mathbb{Z}\left[1 / n\right];
$$
 that is, $(n-1)$ divides  $k$,  which is a contradiction. Hence, there exists  an isomorphism 
$$
(n-1)^i\mathbb{Z}\left[1 / n\right] /(n-1)^{i+1} \mathbb{Z}\left[1 / n\right]
\cong \mathbb{Z}_{n-1}.
$$
The theorem is proved.

\medskip

From this theorem follows the corollary which was proved in  \cite{Mol4}.

 COROLLARY. 1) {\it The group  $BS(1, n)$ is residually nilpotent if and only if 
$n \not= 2$.}

2) {\it The group  $BS(1, 2)$ is not residually nilpotent and the length of its lower central series is equal to   2, i.~e.
$$
\gamma_2 BS(1, 2) = \gamma_3 BS(1, 2).
$$
}

PROOF.
 Suppose that $n \not= 2$. It is enough to prove that  
$\bigcap\limits_{i \geq 2}(n-1)^{i-2}\mathbb{Z}[1 / n]=0$.
Let $\alpha \in \bigcap\limits_{i \geq 2}(n-1)^{i-2}\mathbb{Z}[1 / n]$.
Then $\frac{\alpha}{(n-1)^{i-2}} \in \mathbb{Z}[1 / n]$
for all  $i \geq 2$, that is impossible, since the denominator of the fraction  $\frac{\alpha}{(n-1)^{i-2}}$
  increases  without bound with respect to  as $i$ increases.

 For $n=2$ we have
$$
G=BS(1,2)=\left\langle \, t, a \, \| \, t^{-1}at=a^2 \, \right\rangle .
$$
Since $a=[a,t]$, then $a\in \gamma_s G$, $s\geq 2$. Hence, $a\in \gamma_\omega G$.
Further, 
$$
\left\langle \, t, a \, \| \, t^{-1}at=a^2,\,\, a=1 \, \right\rangle \cong \mathbb{Z}
$$
and $\gamma_2 G\leq \gamma_\omega G$.
Hence
$$
\gamma_2 G = \gamma_\omega G,
$$
and, since  $\gamma_2 G \not= 1$, the group  $BS(1,2)$ is not residually nilpotent.
The corollary is proved.

\medskip

{\bf Non-solvable  Baumslag-Solitar groups.}
The non-solvable Baumslag-Solitar groups have the form 
 $BS(m, \varepsilon m)$, $m > 1$, $\varepsilon = \pm 1$.
In \cite{Mol4}  it was proved that   for $m > 1$, $\varepsilon = \pm 1$  $BS(m, \varepsilon m)$
is residually nilpotent if and only if $m$ is a power of some prime number. 

\medskip

The following proposition describes the lower central series for the Baumslag-Solitar groups  $BS(m, \varepsilon m)$, $m > 1$, $\varepsilon = \pm 1$.

\medskip

PROPOSITION 2.
{\it Let  $m > 1$ be a natural number.
The following isomorphisms hold:

1) $\gamma_i \left(BS(m, m) \right) \cong \gamma_i (\mathbb{Z} * \mathbb{Z}_m),~~i > 1$;

2) $\gamma_i \left(BS(m, -m) \right)\cong \langle a^{2^{i-1}m}\rangle
\times\gamma_i (\mathbb{Z} * \mathbb{Z}_m),~~i > 1.$
}

PROOF.
At first, consider the group  $G=BS(m, m)$.

The natural homomorphisms 
$$
\varphi_i:\gamma_i G \longrightarrow \gamma_i (\mathbb{Z} * \mathbb{Z}_m),~~i \geq 1,
$$
are induced by the  homomorphism 
$$
\varphi : G \longrightarrow \mathbb{Z} * \mathbb{Z}_m
$$
which sends the element  $a^m$ to the unit element in in $\mathbb{Z} * \mathbb{Z}_m$) . Since $a^m$ lies in the center,  
$\mathrm{Ker}\, \varphi_2$ is equal to the subgroup  generated by  $a^{m}$,
and we have the short exact sequence 
$$
1 \to \langle a^m \rangle \to G \to \mathbb{Z} * \mathbb{Z}_m \to 1.
$$
To prove that  $\varphi_i$, $i\geq 2$ is an isomorphism, it is enough to prove that 
$\mathrm{Ker}\, \varphi_2 =1$, i.~e.
$ \gamma_2 G \bigcap \langle a^{m}\rangle =1$.

Since the commutator subgroup  $\gamma_2 G$ of $G$ is  the  homomorphic image of the commutator subgroup of the  free group of rank 2, it  is generated by the commutators 
$$
[t^k,a^l], \quad k,l \in \mathbb{Z}\backslash \{ 0 \}.
$$
Using the fact that  $a^{m}$ lies in the center of  $G$, we can assume that 
$1\leq l \leq m-1$.

The commutator subgroup  $\gamma_2 (\mathbb{Z} * \mathbb{Z}_m)$ of the free product 
$$
\mathbb{Z} * \mathbb{Z}_m=\left\langle \, t, a \, \| \, a^m=1 \, \right\rangle
$$
is the free group with free generators 
$$
[t^k,a^l], \quad k \in \mathbb{Z}\backslash \{ 0 \}, \,\, 1\leq l \leq m-1.
$$
 It follows that  there exists  a natural homomorphism 
$$
\psi : \gamma_2 (\mathbb{Z} * \mathbb{Z}_m) \longrightarrow \gamma_2 G,
$$
 inverse to  $\varphi_2$.
Hence, $\varphi_2$ is an isomorphism  and all
$\varphi_i$, $i\geq 2$, are isomorphisms.

Consider the group  $G=BS(m, -m)$.
 Let us show that there exists an isomorphism 
$$
\gamma_2 G\cong \langle a^{2m}\rangle \times \gamma_2 (\mathbb{Z} * \mathbb{Z}_m).
$$
The commutator subgroup of  $G$ is generated by the commutators
$$
[t^k,a^l], \quad k,l \in \mathbb{Z}\backslash \{ 0 \}.
$$
It is easy to check that for all  $k,l \in \mathbb{Z}$ the following equality holds
$$
[[t^k,a^l],a^m]=1.
$$
Using the commutator identity 
$$
[z,xy]=[z,y][z,x]^y,
$$
we  see that the commutator subgroup of  $G$ is generated by the commutators 
$$
[t^k,a^l], \quad k \in \mathbb{Z}\backslash \{ 0 \}, \,\, 1\leq l \leq m-1
$$
and  powers of the word $a^{2m}$. As above, we conclude that these commutators generate a free group which is isomorphic to the commutator subgroup of the free product 
$\mathbb{Z} * \mathbb{Z}_m$. Hence, the homomorphism 
$$
\varphi_2: \gamma_2 G \longrightarrow \gamma_2 (\mathbb{Z} * \mathbb{Z}_m)
$$
 is invertible. Using the fact that  $a^{2m}$ lies in the center of  $\gamma_2 G$, we  obtain the isomorphism 
$$
\gamma_2 G\cong \langle a^{2m}\rangle \times\gamma_2 (\mathbb{Z} * \mathbb{Z}_m).
$$
Note also the fact that the group  $\gamma_2 (\mathbb{Z} * \mathbb{Z}_m)$
is embedded into the group  $\gamma_2 G$, i.~e. instead of the  isomorphism  symbol we  may write  that of equality.

Suppose that for some 
 $i\geq 2$ we have proved  the existence of the isomorphism 
$$
\gamma_i G\cong \langle a^{2^{i-1}m}\rangle \times\gamma_i (\mathbb{Z} * \mathbb{Z}_m)
$$
and the group  $\gamma_i (\mathbb{Z} * \mathbb{Z}_m)$
is embedded into the group  $\gamma_i G$.
Then
$$
\gamma_{i+1} G =[\gamma_i G, G]=[\langle a^{2^{i-1}m}\rangle \times\gamma_i (\mathbb{Z} * \mathbb{Z}_m), G]=
[\langle a^{2^{i-1}m}\rangle,G] \times [\gamma_i (\mathbb{Z} * \mathbb{Z}_m), G]=
$$
$$
=\langle a^{2^{i}m}\rangle \times [\gamma_i (\mathbb{Z} * \mathbb{Z}_m), \mathbb{Z} * \mathbb{Z}_m]=
\langle a^{2^{i}m}\rangle \times\gamma_{i+1} (\mathbb{Z} * \mathbb{Z}_m).
$$
The proposition is proved.

\medskip


Thus, Theorem 1 and Proposition  2 imply  the following criteria of residual nilpotence of Baumslag-Solitar groups,  found in   \cite{Mol4}.

\medskip

THEOREM 2.
{\it The group  $BS(m,n)$, $0 < m \leq |n|$, is residually nilpotent if and only if  $m = 1$ and $n \not= 2$, or $|n| = m > 1$ and $m$ is a  power of some prime number.}

\medskip

Denote by  $r \mathfrak{F} BS$ the set of RF 
 Baumslag-Solitar groups, 
by $r \mathfrak{N} BS$, the set of RN  Baumslag-Solitar groups,
by $r \mathfrak{F}_p BS$, the set of RF$_p$ Baumslag-Solitar groups,
by $\{r \mathfrak{F}_p BS\}_p$, the union $r \mathfrak{F}_p BS$ by all prime  $p$.
The following inclusions hold
$$
r \mathfrak{F}_p BS \subset \{r \mathfrak{F}_p BS\}_p \subset r \mathfrak{N} BS \subset r \mathfrak{F} BS.
$$
The difference  $ r \mathfrak{F} BS \setminus r \mathfrak{N} BS$ contains 
the group $BS(1, 2)$ and groups $BS(m, \pm m)$, where $m$ is not a power of some prime number.
The difference $ r \mathfrak{N} BS \setminus \{r \mathfrak{F}_p BS\}_p$
contains the set of groups  $BS(p^r, -p^r)$, where $p > 2$ is a prime number.

\bigskip

\begin{center}
{\bf \S~3. Intersection of the terms of the lower central series }
\end{center}

\bigskip

As we know, if a group $G$ is residually nilpotent, then  $\gamma_{\omega} G = 1$. In the opposite case, the quotient   $G / \gamma_{\omega}$
is the maximal quotient of $G$ which is residually nilpotent.
In this section we find  intersections of all terms of the lower central series for  several Baumslag-Solitar groups.

The proof of the following lemma is evident. 

\medskip

LEMMA 1. {\it The abelianization  $BS(m,n)^{ab}$ is isomorphic to $\langle a, t ~||~a^{n-m} = 1\rangle$. In particular:

-- if $n = m+1$, then $BS(m,n)^{ab} \cong \mathbb{Z}$;

-- if $n = m$, then $BS(m,n)^{ab} \cong \mathbb{Z} \times \mathbb{Z}$;

-- if $n \not= m, m+1$, then $BS(m,n)^{ab} \cong \mathbb{Z} \times \mathbb{Z}_{n-m}$.
}

\medskip

We also need

\medskip

LEMMA 2. {\it In the group $BS(m,n)$ for all $i \geq 0$ the following inclusion holds  $a^{(n-m)^i} \in \gamma_{i+1} BS(m,n)$.}

PROOF. From the defining relation of  $BS(m,n)$ follows the relation 
$$
[a^m, t] = a^{n-m},
$$
i.~e. $a^{n-m} \in \gamma_2 BS(m,n)$. Raising both parts of the defining relation 
of $BS(m,n)$  to the power  $n-m$, we get
$$
t^{-1} a^{m(n-m)} t = a^{n(n-m)} \Leftrightarrow a^{m(m-n)} t^{-1} a^{m(n-m)} t = a^{(n-m)^2}.
$$
Hence,
$$
[[a^m, t]^n, t] = a^{(n-m)^2}.
$$
Using induction  on  $i$, we get the  sought  inclusion.  The lemma is proved.

\medskip

Further, we shall consider the Baumslag-Solitar groups which are not residually nilpotent.

 Amongst non-abelian Baumslag-Solitar groups the group  $BS(m,m+1)$, $m>1$,
has the  shortest lower central series.  

\medskip

PROPOSITION 3. {\it The  lower central series of the group  
$$
G=BS(m,m+1)=\left\langle \, a, t \, \| \, t^{-1} a^m t=a^{m+1} \, \right\rangle,~m > 1,
$$
 has length 2.}

PROOF.
Since $a = [a^m,t]$, then $a\in \gamma_s G$, $s\geq 2$.  It follows that $a\in \gamma_\omega G$.
 Moreover,
$$
\left\langle \, a, t \, \| \, t^{-1} a^m t = a^{m+1},\,\, a=1 \, \right\rangle \cong \mathbb{Z}.
$$
 Hence $\gamma_2 G\leq \gamma_\omega G$ and 
 $\gamma_2 G = \gamma_\omega G.$
The proposition is proved.

\medskip

 Setting $d =$GCD$(m, n)$,  the greatest common devisor  of $m$ and $n$, 
 let us consider the group
$G=BS(m,n)$ for $m=kd$, $n=kd + d$, where $k$, $d \in \mathbb{Z}$.
In $G$  the following relation holds
$$
a^{d}=[a^{kd}, t]\in \gamma_2 G.
$$
Hence, $a^d\in \gamma_\omega G$.
The next theorem generalises Proposition 3.

\medskip

THEOREM 3.
{\it In $G=BS(kd, kd + d)$, $k$, $d \in \mathbb{Z}$ the following statements are true: 

1) if $d$ is a power of some prime number or  $d=1$, then
$\gamma_\omega G =\left\langle \, a^{d} \, \right\rangle^{G};$

2) if $d$ is not a power of some prime number, then
$\gamma_\omega G > \left\langle \, a^{d} \, \right\rangle^{G}.$
}

PROOF.
There exists an  isomorphism
$$
G/ \langle \, a^d \, \rangle^{G}=
\langle \, a, t \, \| \, a^d=1 \, \rangle \cong \mathbb{Z}\ast \mathbb{Z}_d.
$$
If $d$ is a power of some prime number or  $d=1$, then
the free product  $\mathbb{Z}\ast \mathbb{Z}_d$ is residually nilpotent and we have 
$\gamma_\omega (G/\langle \, a^d \, \rangle^{G})=1$
Hence, the following inclusion holds 
$$
\gamma_\omega G \subseteq \langle \, a^d \, \rangle^{G}.
$$
The inverse inclusion was proved before. 

If $d$ is not a power of  a  prime number, then
$\gamma_\omega (\mathbb{Z}\ast \mathbb{Z}_d)\neq 1$ and
$\gamma_\omega G > \left\langle \, a^{d} \, \right\rangle^{G}.$

The theorem is proved.

\medskip

Let us fix the Baumslag-Solitar group  $G = BS(m, n)$, $d =$GCD$(m, n)$, which is not residually finite.
Let $R$ be the intersections of all subgroups of finite index in  $G$. D.~I.~Moldavanskii \cite{Mol2} proved that
$R$ is the normal closure of the set of commutators 
$$
[t^k a^d t^{-k}, a],~~k \in \mathbb{Z},
$$
in $G$. Let  $\overline{G} = G / R$ be the maximal residually finite quotient of  $G$. Then
$$
\overline{G} = \langle a, t ~|~t^{-1} a^m t = a^n,~~[t^k a^d t^{-k}, a]=1,~~k \in \mathbb{Z} \rangle.
$$
Define  the subgroup $\overline{A}\subset\overline{G}$,
$$
\overline{A} = \langle t^k a^d t^{-k}, ~~k \in \mathbb{Z} \rangle.
$$
Then $\overline{A}$ is  an abelian normal subgroup of  $\overline{G}$ and
we have the isomorphism  (see \cite{Mol2})
$$
\overline{G} / \overline{A} = \langle a, t~|~a^d = 1 \rangle = \mathbb{Z}_d * \mathbb{Z}.
$$

It is interesting to find  a connection between  $R$ and $\gamma_\omega G$.

\medskip

Note that if  $d \geq 2$, then $BS(m, n)$ contains a free non-abelian subgroup. More precisely, we can prove

\medskip

PROPOSITION 4. {\it If $d \geq 2$, then $G = BS(m, n)$ contains a free non-abelian subgroup  of infinite rank.}

\medskip

PROOF.
Let $D$ be the Cartesian subgroup of  $\mathbb{Z}_d * \mathbb{Z}$; that is, the kernel of the natural homomorphism  $\mathbb{Z}_d * \mathbb{Z} \longrightarrow \mathbb{Z}_d \times \mathbb{Z}$.
It is clear that $D=(\overline{G} / \overline{A})'$.
 It is well known that  $D$ is the free group with free generators
$$
[t^k,a^s], \quad k \in \mathbb{Z}, \quad 1 \leq s \leq d-1.
$$
Hence, if $d\geq 2$, then $\overline{G}$, and, correspondingly,  $G$, contains a free non-abelian subgroup of infinite rank. The proposition is proved. 

\medskip

Note that $R \subseteq G'$ and Lemma 1 implies that 
$$
G/G' \cong \mathbb{Z} \times \mathbb{Z}_{n-m}.
$$
On the other  hand,
$$
\overline{G}/\overline{G}'\overline{A} \cong \mathbb{Z} \times \mathbb{Z}_{d}.
$$
If $n-m > d$, it implies that 
$\overline{G}'$ is a proper subgroup of  $\overline{G}'\overline{A}$.
If $n-m = d$, then $\overline{G}'=\overline{G}'\overline{A}$, i.~e. we have the inclusion 
$\overline{A} \subseteq \overline{G}'$.

Since the quotient  $G/\gamma_\omega G$ is residually nilpotent, the intersection of all its subgroups of finite index is trivial.  Hence, $R$
lies in  $\gamma_\omega G$ and we have the inclusion
$$
R=\left\langle [t^k a^d t^{-k}, a],~~k \in \mathbb{Z}\right\rangle^G \subseteq \gamma_\omega G.
$$
Let us denote  $\gamma_\omega G / R= \overline{\gamma_\omega G}.$
Note that $\overline{\gamma_\omega G} < \gamma_\omega\overline{ G}$.

Let $A \leq G$ be a such subgroup that  $A/R=\overline{A}$. We have  the inclusions 
$$
G \geq G'A \geq A \geq R,
\quad
G \geq G' \geq \gamma_\omega G \geq R
$$
with  quotients 
$$
G / G'A \cong \mathbb{Z} \times \mathbb{Z}_{d},
\quad
G / A \cong \mathbb{Z} \ast \mathbb{Z}_{d},
\quad
G / G' \cong \mathbb{Z} \times \mathbb{Z}_{n-m},
$$
$$
G'A/ A \cong
\left\{\begin{array}{cc}
F_{\infty}, & \mbox{if} \,\, d \geq 2,\\
\mathbb{Z}, & \mbox{if} \,\, d =1. \\
\end{array}
\right.
$$
Also note that the quotient  $\gamma_\omega G / R $ is abelian.

If $n=m+d$, i. e. $m=kd,$ $n=kd+d$, then as we show above 
 $A \subseteq G'$ and
$$
G \geq G' \geq A \geq R.
$$
In particular, $\gamma_2 G / \gamma_\omega G$ contains a free abelian subgroup of infinite rank  for $d \geq 2$.

If $d$  is some power of a prime number, then the quotient   $\overline{G}/\overline{A} \cong \mathbb{Z} \ast \mathbb{Z}_{d}$
is residually nilpotent, $\gamma_\omega\overline{ G} \subseteq \overline{A}$ and the following inclusions hold
$$
G \geq G'A \geq A \geq \gamma_\omega G \geq R.
$$

If $n=m+d$ and $d$  is a power of some prime number, then 
$$
G \geq G' \geq A \geq \gamma_\omega G \geq R.
$$

If $n=m+1$, then the Proposition 3 implies that  $G' = A = \gamma_\omega G$.

Hence, we  have proved

\medskip

PROPOSITION 5.
{\it
The following inclusions hold
$$
G \geq G'A \geq A \geq R,
\quad
G \geq G' \geq \gamma_\omega G \geq R,
$$
where the quotient  $A / R $ is abelian. In  particular,

1) If $n-m\geq d>1$, then
$$
G / G'A \cong \mathbb{Z} \times \mathbb{Z}_{d},
\quad
G / A \cong \mathbb{Z} \ast \mathbb{Z}_{d},
\quad
G / G' \cong \mathbb{Z} \times \mathbb{Z}_{n-m},
\quad
G'A/ A \cong F_{\infty}.
$$

2) If $n-m> d=1$, then $A$ contains $G'$ as a proper subgroup and 
$$
G / A \cong \mathbb{Z},
\quad
A / G' \cong \mathbb{Z}_{n-m}.
$$

3) If $n-m= d=1$, then $G' = A = \gamma_\omega G$.

4) If $d$  is some power of a prime number, then
$$
G \geq G'A \geq A \geq \gamma_\omega G \geq R.
$$

5) If $n=m+d$ and  $d$ is a power of some prime number, then 
$$
G \geq G' \geq A \geq \gamma_\omega G \geq R.
$$

}

\medskip

REMARK.
If $n-m= d=1$, then $\gamma_\omega G / R \cong \mathbb{Z}^{\infty}$
and
$$
\overline{G}\cong \mathbb{Z} \wr \mathbb{Z}
$$
is the wreath product of two infinite cyclic groups and this group  is isomorphic to the group that is generated by the matrices 
$$
\left(%
\begin{array}{cc}
1 & 1 \\
0 & 1 \\
\end{array}%
\right),\quad
\left(%
\begin{array}{cc}
\xi & 0 \\
0 & 0 \\
\end{array}%
\right),
$$
where  $\xi$ is a  transcendental number.

\medskip

Proposition 5 raises the  questions.

{\bf Question 1.}  Can we find  the structure of the quotient
$\gamma_\omega G/R$ for $n-m\geq d>1$ and $n-m> d=1$?

{\bf Question 2.} Find   connections between the subgroup 
$A$ and $\gamma_\omega G$, if $d$ is not a power of some prime number. 

\medskip

The following questions also seem interesting.

{\bf Question 3.} Is it true that if $G$ is a Baumslag-Solitar group,
then $\gamma_{\omega} G = \gamma_{\omega+1} G$, where $\gamma_{\omega+1} G = [\gamma_{\omega} G, G]$?

Recall  (see, for example \cite{F}) that a finitely generated group, acting on a tree  in such a way that all stabilizers of  vertices and edges are infinite cyclic groups 
 is called a  {\it generalized Baumslag-Solitar group}.

{\bf Question 4.} Which generalized Baumslag-Solitar groups are residually nilpotent?


\vspace{1cm}

\begin{thebibliography}{HD}


\bibitem{Ba}
V.  Bardakov, On the width of verbal subgroups of some free constructions (Russian), Algebra i
Logika 36, (1997),  494--517; translation in Algebra and Logic 36,  (1997),  288–-301.


\bibitem{BS}
G. Baumslag, D. Solitar, Some two-generator one-relator non-Hopfian groups,
Bull. Amer. Math. Soc., 68,  1962, 199--201.

\bibitem{F}
M. Forester, Splittings of generalized Baumslag-Solitar groups,
Geometriae Dedicata, 121, 1 (2006), 43--59.

\bibitem{Mes}
S. Meskin, Nonresidually finite one-relator groups, Trans. Amer. Math. Soc., 164,
1972,  105--114.

\bibitem{Mal}
 A. I. Malcev, On homomorphisms onto finite groups, Uchen. Zapiski Ivanovsk.
ped. instituta, 18, 5 (1958), 49–60 (also in “Selected papers”, Vol. 1, Algebra,
1976, 450--462) (Russian).

\bibitem{M}   A. I. Malcev, Generalized nilpotent algebras and their associated groups,
Mat. Sbornik N.S., 25, 3 (1949), 347-–366.

\bibitem{Mol1}
D. I. Moldavanskii, On $p$-residually finiteness  of  HNN-extensions. (Russian) Vestnik Ivanov. Gos. Univ.,  3 (2000), 129--140.

\bibitem{Mol2}
D. I. Moldavanskii, On the intersection of subgroups of finite index in the Baumslag-Solitar groups. (Russian) Mat. Zametki 87 (2010),  1, 92-–100; 
translation in Math. Notes,  87,  1-2 (2010),  88-–95.


\bibitem{Mol4}
D. I. Moldavanskii,  Residual Nilpotence of Groups with One Defining Relation. (Russian) Mat. Zametki, 107, 5 (2020),  752-–759.


\bibitem{Mol}
D. I. Moldavanskii, On some residual properties of Baumslag-Solitar groups, ArXiv, 2013.


\bibitem{Mol3}
D. I. Moldavanskii, On the residual properties of Baumslag-Solitar
groups, Communications in Algebra, 46,  9 (2018), 3766--3778.




\end{thebibliography}
\end{document}